\documentclass[12pt]{article}
\title{Factored Notation for Interval I/O
\thanks{Report DCS-264-IR, Department of Computer Science,
University of Victoria, Victoria, BC, Canada.}
}
\author{\mbox{M.H. van Emden}\\
        \mbox{Aerospace Faculty, TU Delft and 
              Constraints Group, CWI}\\
        \mbox{On leave from Computer Science Dept,
              University of Victoria}
	}
\date{February 15, 2001}
\begin{document}
\maketitle

\begin{abstract}
This note addresses the input and output of intervals in the sense of
interval arithmetic and interval constraints.  The most obvious, and
so far most widely used notation, for intervals has drawbacks that we
remedy with a new notation that we propose to call {\em factored
notation}.  It is more compact and allows one to find a good trade-off
between interval width and ease of reading.  We describe how such a
trade-off can be based on the information yield (in the sense of
information theory) of the last decimal shown.
\end{abstract}

\section{Introduction}

Once upon a time, it was a matter of professional ethics among computers
never to write a meaningless decimal. Since then computers have become
machines and thereby lost any form of ethics, professional or
otherwise.  The human computers of yore were helped in their ethical
behaviour by the fact that it took effort to write spurious decimals.
Now the situation is reversed: the lazy way is to use the default
precision of the I/O library function. As a result it is common to see
fifteen decimals, all but three of which are meaningless.

Of course interval arithmetic is not guilty of such negligence. After
all, the very {\em raison d'\^etre} of the subject is to be explicit
about the precision of computed results.  Yet, even interval
arithmetic is plagued by phoney decimals, albeit in a more subtle
way. Just as conventional computation often needs more care in the
presentation of computational results, the most obvious interval
notation with default precision needs improvement.

As a bounded interval has two bounds, say, $l$ and $u$, the most
straightforward notation is something like $[l,u]$.  Written like
this, it may not be immediately obvious what is wrong with writing it
that way.  But when confronted with a real-life consequence\footnote{
As can be found on page 122 of \cite{vhlmyd97}.  }, such as the
statement that an unknown real $x$ belongs to
\begin{eqnarray}
      [+0.6180339887498946804,+0.6180339887498950136], \label{crude}
\end{eqnarray}
it becomes clear that this is not a practical notation.

This problem has been primarily addressed by Schulte, Zelov, Walster,
and Chiriaev \cite{szwc99}, who refer to the $[l,u]$ notation as the
``conventional one'', and observe:
\begin{quote}
With conventional interval output, it is often difficult to determine
the relative sharpness of interval results, $\ldots$ This is especially
true when large amounts of interval data need to be examined.
\end{quote}

\subsection{The scanning problem}
To start with, the notation exemplified by (\ref{crude}) requires one
to scan both bounds to find the leftmost different digit.  Only then
does one have an idea of the width of the interval.  In casting about
for alternative notations, we note that an interval has more
potentially useful attributes than the lower bound $l$ and the upper
bound $u$.  We could use the centre $c$, the radius $r$, or the
diameter $d$.  Any two independent attributes could be paired to give
a notation for intervals.  Of course, several such combinations are
awkward. A reasonable alternative to pairing $l$ and $u$ is to
combine $c$ with $r$.  But then, in writing the pairs, we should drop
the suggestive ``['' and ``]''.  Thus we get $\langle c,r \rangle$. In
this notation, for example, the interval~(\ref{crude}) becomes:
\begin{eqnarray}
\langle 0.6180339887498948470,0.0000000000000001666 \rangle.
                                \label{centered}
\end{eqnarray}
The advantages of the $\langle c,r \rangle$ notation include solving
the problem of having to scan for the leftmost differing digit. Also,
it shows at a glance the width of the interval.  While this is an
improvement over (\ref{crude}), it is still awkward that instead of a
long string of repeated digits, one now has a long string of zeros.

\paragraph{Factored notation as solution to the scanning problem}
Now the problem with (\ref{crude}) is repeated occurrence of
characters, much like the repeated occurrence of $a$ in the expression
$ab + ac$. This analogous problem is solved by {\em factoring}, so
that one gets $a(b+c)$. Consider now string concatenation as an
infix multiplicative operation. That suggests ``factoring'' in
(\ref{crude}) the longest initial sequence of common digits, so that
$$ x \in [0.6180339887498946804,0.6180339887498950136] $$ becomes
\begin{eqnarray}
         x \in 0.61803398874989[46804,50136] \label{factored}
\end{eqnarray}
We call this the {\em factored notation}.
It solves the scanning problem in a more compact way than the
$\langle c,r \rangle$ notation.

\subsection{The problem of useless digits}
By solving the scanning problem, factored notation brings out the
other problem: does one, in (\ref{factored}) for example, really want
to know the width of an interval in a relative precision of $10^{-5}$?
One should keep in mind that relaxing the relative precision in the
{\em width} does not detract from the $10^{-14}$ relative precision in
$x$ itself.  If it is supremely important to preserve {\em all} the
information that is available about $x$, then one should indeed
preserve all the digits inside the brackets.  But usually one does not
want to keep all information: what matters is how much information the
human reader can absorb in whatever limited time is available.  By
giving too much information about the width of the interval, it is
usually the case that {\em less} information about $x$ is actually
transferred.  So one should present an {\em approximation} to the
available information about $x$ that can be represented with fewer
digits.

Now approximation is an ambiguous concept.  Here we mean it in the
sense of giving less information, under the constraint that what
information remains truly is information, that is, a true statement.
In this paper we assume that the purpose of using intervals at all is
to obtain from a computation a true statement about the real-valued
components of a solution.  We assume that these statements say that a
certain real belongs to an interval and that this interval is
interpreted as a set of reals.  However wide the interval $I$ is,
there is nothing approximate about the statement $x \in I$ as such.
Each bound, though a mere machine number, is also an unambiguously
defined real. So $I$ is an unambiguously defined set.  The relation
$\in$ either holds or does not hold.

Computational results of the form $x \in I$ are made possible by
correct rounding and by methods such as Interval Newton
\cite{hnsn92} or Interval Constraints \cite{vhlmyd97}
that take into account errors from all sources.  Suppose that the $I$
in $x \in I$ is the raw computational result: an interval with binary
floating-point numbers as bounds. For output, we need to represent it
as a pair of decimal numerals in factored notation. In general,
rounding is necessary in the transition between binary and decimal and
care has to be taken to ensure that the rounding is outward.

Of course we want to keep the resulting interval as small as possible:
the smaller, the more information about $x$.  However, it may be that
in the chosen output language for the interval (decimal numerals), the
description of $I$ requires many numerals, in factored notation.  Then
we may prefer instead an output $x
\in I^\prime$, with $I \subseteq I^\prime$, such that $I^\prime$ has a
simple representation. In that sense $I^\prime$ is an approximation to
$I$.  But there is nothing approximate about the truth value of $I
\subseteq I^\prime$, and hence of $x \in I^\prime$.

This suggests a sequence of intervals
\begin{eqnarray}
I \subseteq I_1 \subseteq I_2 \subseteq I_3 \ldots \label{seq}
\end{eqnarray}
where each following interval requires fewer decimal digits than the
one before.

As an illustration consider the interval that, in $[l,u]$ notation
would appear as $$ [5.1268427635136,5.1268472635136].  $$

In the example the sequence $I, I_1, I_2, I_3, \ldots$, in factored
notation, is as shown in Figure~\ref{pyramid}.

\begin{figure} \label{pyramid}
\begin{eqnarray*}
5.12684[27635136,72635136]     & \subseteq &  \\
 5.12684[2763513,7263514]      & \subseteq &  \\
  5.12684[276351,726352]       & \subseteq &  \\
   5.12684[27635,72636]        & \subseteq &  \\
    5.12684[2763,7264]         & \subseteq &  \\
     5.12684[276,727]          & \subseteq &  \\
      5.12684[27,73]           & \subseteq &  \\
       5.12684[2,8]            & \subseteq &
\end{eqnarray*}
\caption{Successive approximations to an interval. Which is best?}
\end{figure}

\section{Semantics of decimal numerals} \label{semDec}
We need decimal numerals to denote the bounds in an interval. Now,
what exactly is meant by a given numeral?  There seem to be two
interpretations.

The {\em first} interpretation is appropriate in the context of
physical measurements, where the convention is used that, in the
absence of explicit qualification, all digits are significant. For
example, when $y = 0.123$ is given as outcome of an experiment, the
implication is that the fourth and later decimals may have any value.
That is, no more is implied than $0.1230\infty \leq y \leq
0.1239\infty$ where $\infty$ indicates infinite repetition of the
preceding decimal.

This first interpretation has drastic consequences when we consider
$y^2$, which cannot be given as $0.015129$, because the last three
digits are likely to be incorrect as the next digit after the ones
given in $0.123$ is not known.

One way of expressing this interpretation of $y=0.123$ is to say that
$y$ belongs to a certain interval. It is natural to express this
interval as $[a,b]$. But how to write $a$ and $b$?  If we write
$[0.123,0.124]$, then these numerals need to be interpretated in a
different way.

Hence the {\em second} interpretation of a decimal numeral: the one that
implies infinitely many zeros after the last decimal.

Thus, $x=0.d_1 \ldots d_k$ means either
$$ 0.d_1 \ldots d_k0\infty \leq x \leq 0.d_1 \ldots d_k9\infty $$
or
$$ x = 0.d_1 \ldots d_k0\infty $$

The first interpretation we call {\em physics convention}. It can be
regarded as a restricted interval notation; one where only intervals
can be denoted in which the bounds differ by one unit in the last
decimal.

To implement our chosen semantics of intervals in general, we need to
write decimals in the bounds. But for these the physics convention is
not satisfactory; we need the latter interpretation: bounds of
intervals, in so far as they exist, are reals, not sets of reals. We
call this the {\em point convention}.

\section{How many digits for interval bounds?} \label{numDigits}
According to Shannon's theory of information, observations can reduce
the amount of uncertainty about the value of an unknown quantity. The
amount of information yielded by an observation is the decrease (if
any) in the amount of uncertainty. Shannon argues that the amount of
uncertainty is appropriately measured by the entropy of the
probability distribution over the possible values. For a uniform
distribution on a finite number of values, this reduces to the
logarithm of the number of possible values.  For an unknown real $x$
that belongs to an interval we assume a uniform distribution over the
interval and we take for the uncertainty the logarithm of the width of
the interval.

The amount of information gained by resolving the choice between two
equally probable alternatives is called the $b$inary u{\it it} of
information, or, the {\it bit}.
This is widely used. Strictly analogously, but less widely used, is the
$d$ecimal un{\it it} of information, or, the {\it dit}.

Before we look at intervals, let us evaluate, from an
information-theoretic point of view, numerals interpreted according to
the physics convention. There $x = 0.123$ means $x \in [0.123,0.124]$
with $\log _{10}0.001 = -3$ as measure of uncertainty. Add one
decimal, and the measure of uncertainty is reduced to $-4$. The
information yielded is one dit.  At the same time, $0.123$ leaves as
uncertainty which the next more precise numeral is: $$0.1230\hbox{ or
} 0.1231 \hbox{ or } \ldots
\hbox{ or } 0.1239.$$ These are ten possibilities, which are reduced
to one by supplying the fourth decimal. Thus we see a perfect match
between the expense (one dit) of adding a digit and the yield of
information (one dit). Such a good match is rarely the case with
interval notation.

In (\ref{seq}) we described how one can save on decimals while
maintaining inclusion by widening the interval.  To judge whether it
is worthwhile to display all available decimals, let us determine the
amount of information contained in the last decimal in each of the
bounds for an interval.

Suppose all available decimals are as given in $$ x \in
5.12684[27635136,72635136].$$ Let us consider the amount of
information about $x$ that is communicated by the successive digits
within the brackets, starting with the second digit. If we only had
one digit between the brackets, then we would write $ x \in
5.12684[2,8]$, with an uncertainty in $x$ equal to $\log _{10} 6
\times 10^{-6} = -6 + \log _{10} 6$.  The next digit gives $ x \in
5.12684[27,73]$, with an uncertainty of $\log _{10} 46 \times 10^{-7}
= -6 + \log _{10} 4.6$.  The decrease in uncertainty, that is the
information, given by the last digit inside the brackets is $\log_{10}
6 - \log_{10} 4.6$, which is considerably less than $1$.

Jumping to the final digit now: without it we have for the uncertainty in
$$ x \in 5.12684[27635136,72635136]$$
an amount in dits equal to
$\log_{10} 4.5000000  \times  10^{-6}$.
If we had to do without that last digit, then all we can say is
$$ x \in 5.12684[2763513,7263514]$$
with an uncertainty equal to 
$- \log_{10} 4.500001  \times  10^{-6}$.

Using the approximation that $\log(1+\alpha)$ is in the order of
$\alpha$ for small $\alpha$, we find that the difference is in the
order of $10^{-6}$ dits. This is the yield of writing the last digit
in each bound. Ideally, the yield is one dit. With the physics
convention for a decimal numeral, this maximum possible yield is
guaranteed.  As soon as we use the more flexible interval notation, we
have to accept a lower yield. But accepting a yield as low as
$10^{-6}$ dits is most situations far from optimal.

If one insists on the maximum yield possible with interval notation,
one would write $x \in 5.12684[2,8]$, which is a perfectly good
summary of what is known about $x$. One may have sympathy for those
who point out that the omitted decimals represent information that was
gained at some, possibly considerable, computational cost. In that
case, one would add one or two decimals, as in $ x \in
5.12684[276,727]$. The next decimal will only add a negligible amount
of information. Even though $ x \in 5.12684[27635136,72635136]$ might
represent the investment of huge computational resources, they would
be wasted, as there just is hardly any more information than in the
three-decimal version.

This is then our recommendation: use the factored notation to solve the
scanning problem, and as default give two or three decimals inside the
brackets.

\section{The general case}

In general one may have a bound written in scientific notation; that
is, with an explicit power of ten. In case both bounds are so written,
we propose to normalize one mantissa and to factor out its power of
ten from both bounds. In the remaining bound, the decimal point is
then adjusted so that its power of ten is zero, and can therefore be
omitted.

For example,
$$ x \in [5.1268427635136 \times 10^{2},5.1268472635136 \times 10^{3}] $$
becomes
$$ x \in [0.51268427635136,5.1268472635136] \times 10^{3} $$
Then the above considerations apply. They apply in the strongest way
when both bounds have different normalized powers of ten; then we have
the situation of a wide interval so that only few decimals should be
used inside the brackets; $ x \in [0.51,5.2] \times 10^{3} $ seems
plenty.

Suppose we write an interval as
$$
.d_0 \ldots d_{n-1}[e_0 \ldots e_k, f_0 \ldots f_k]
$$
then we have a width of about $10^{-(n+1)}$, neglecting the factor of
$\log(10)$.  A unit in the last digits $e_k$ and $f_k$ can change the
width by
at most $10^{-(n+1+k)}$.  Hence the width becomes about
$10^{-(n+1)}(1+10^{-k})$.  Thus the uncertainty decreased by the last
digit is about $\log_{10}(1+10^{-k})$, which is about $10^{-k}$. Hence
our opinion that $k$ equal to two or three is plenty.

\section{Related work} \label{relatedWork}
Hansen \cite{hnsn92}, Hammer et al. \cite{hhkr95}, and Kearfott
\cite{krftt96} opt for the straightforward $[l,u]$
notation. Hansen mostly presents bounds with few digits, but for
instance on page 178 we find
$$
	     [0.192895419,0.192895434],
$$
demonstrating the problems addressed here.

The standard notation in the {\em Numerica} book \cite{vhlmyd97} solves
the scanning problem in an interesting way.  It uses the idea of the
$\langle c,r \rangle$ notation, but writes instead $c+[-r,+r]$.  This
variation has the advantage of not introducing new notation.  The
reason why we still prefer factored notation is clear from the
$\langle c,r \rangle$ example (\ref{centered}),
which, if rewritten as $c+[-r,+r]$ becomes
$$
0.6180339887498948470 + [-0.0000000000000001666,+0.0000000000000001666]
$$
Although it is attractive not to introduce special-purpose notation,
there is so much redundancy here that the factored alternative:
$$
0.61803398874989[46804,50136]
$$
seems worth the new notation.

Schulte et al. \cite{szwc99} describe {\em single-number interval I/O}
as remedy for the scanning problem. This notation is used in several
systems listed in \cite{szwc99}. Details differ between input and
output. Here we sketch only the main features shared by single-number
input and output. The three basic formats are $[d]$, $d$, and
$[d_1,d_2]$, where $d$, $d_1$, and $d_2$ are decimal numerals.  $[d]$
stands for the interval $[d,d]$.  Suppose {\em uld} is the real number
that is the value of the unit of the last decimal in $d$, then the
meaning of $d$ in single-number interval notation is $[d - uld, d +
uld]$.  Finally, $[d_1,d_2]$ means just that; it appears that the
tacit assumption is made of infinitely many zeros after the last
digits of $d_1$ and $d_2$.

Single-number interval I/O can only solve the scanning problem by using
the format $d$ (single numeral without brackets).  This implies a
solution to the problem of useless digits, but more drastically than
one might like.  For example, $ 0.123[45678,56789] $ would be
rendered as $0.1235$. Thus the original interval is widened to $
0.123[4,6]. $
One might want to retain more information about the width of the 
interval. In factored notation one can choose among
$0.123[4,6]$,
$0.123[45,57]$,
$0.123[456,568]$,
$\ldots$.

In his interval system {\sc clip} \cite{hck00}, Hickey has a form of
what is called ``single-number interval I/O'' in \cite{szwc99}, as well
as several valuable distinctions not found elsewhere.  Let {\em uld} be
the real number that is the value of the unit of the last decimal in
$d$.  The formats described by Hickey include \verb+d+, \verb+d*+, and
\verb+d...+, where \verb+d+ is a decimal numeral.

\verb+d+ denotes the smallest interval that includes $d$ and has binary
floating numbers as bounds.
\verb+d*+ denotes $[d - uld, d + uld]$.  \verb+d...+ denotes
$[d , d + uld]$.  Finally, in {\sc clip}, \verb+d#+
is used to denote, not an
interval, but a binary floating-point number nearest $d$.  As in
\cite{szwc99}, the scanning problem is solved in {\sc clip}.  But the
problem of useless digits is solved in the same inflexible way as in
\cite{szwc99}.

\section{Conclusions}
Interval methods are coming of age.  When interval software was
experimental, it didn't matter whether interval output was easy to
read.  Now that the main technical challenges have been overcome, and
we at least {\em know} how to ensure that the floating-point bounds
include all reals that are possible values of the variable concerned,
we need to turn our attention to small, mundane matters of
housekeeping, which include taking care of the convenience of users.
Factored notation seems to be an advance in this respect compared to
alternatives described in the literature.

One may not agree with the assumptions made here for evaluating the
utility of the last decimal used in specifying an interval bound.
Information theory needs a probability distribution before it can
assign an amount of information to an observation. The uniform
distribution assumed here is not convincing. But it need not be: the
amount of information in every next decimal declines so steeply, that
most other assumptions would lead one to the same conclusion about a
reasonable number of decimals inside the brackets of factored
notation.

\section{Acknowledgments}
Many thanks to Fr\'ed\'eric Goualard for helpful comments on a draft of
this paper. Support from CWI, NWO, NSERC, TU Delft, and the University
of Victoria is hereby gratefully acknowledged.

\end{document}